\documentclass[12pt]{article}
\usepackage[latin1]{inputenc}
\usepackage{amssymb,amsmath,amsfonts}
\newtheorem{theorem}{Theorem}
\newtheorem{proposition}[theorem]{Proposition}

\newtheorem{corollary}[theorem]{Corollary}

\newcommand{\Q}{\mathbb{Q}}

\newcommand{\Ric}{\mbox{Ric}}

\newcommand{\po}{{\hspace*{-1ex}}{\bf .  }}

\DeclareMathAlphabet{\mathpzc}{OT1}{pzc}{m}{it}

\def\<{\langle}
\def\>{\rangle}

\def\id{I}

\def\bea{\begin{eqnarray*} }
\def\eea{\end{eqnarray*} }
\def\be{\begin{equation} }
\def\ee{\end{equation} }

\def\proof{\noindent{\it Proof:  }}
\def\qed{\ifhmode\unskip\nobreak\fi\ifmmode\ifinner
\else\hskip5 pt \fi\fi\hbox{\hskip5 pt \vrule width4 pt
height6 pt  depth1.5 pt \hskip 1pt }}

\begin{document}

\title{Einstein submanifolds with flat normal bundle in space forms
are holonomic}
\author{M.\ Dajczer, C.-R. Onti and Th.\ Vlachos}
\date{}
\maketitle

\begin{abstract} A well-known result asserts that any isometric 
immersion with flat normal bundle of a Riemannian manifold with 
constant sectional curvature into a space form is (at least locally) 
holonomic.  In this note, we show that this conclusion remains valid 
for the larger class of Einstein manifolds. As an application, when 
assuming that the index of relative nullity of the immersion is a 
positive constant we conclude that the submanifold has the structure 
of a generalized cylinder over a submanifold with flat normal bundle.
\end{abstract}

\renewcommand{\thefootnote}{\fnsymbol{footnote}} 
\renewcommand{\thefootnote}{\arabic{footnote}} 

\renewcommand{\thefootnote}{\fnsymbol{footnote}} 
\renewcommand{\thefootnote}{\arabic{footnote}}

A remarkable class of submanifolds in space forms are those that 
enjoy the property of being  holonomic. An isometric immersion 
$f\colon M^n\to\Q_c^N$ of a Riemannian manifold into a space form of 
constant sectional curvature $c$ is said to be \emph{holonomic} 
if $M^n$ carries a global system of orthogonal coordinates such that 
at any point the coordinate vector fields diagonalize its second 
fundamental form $\alpha\colon TM\times TM\to N_fM$ with values in 
the normal bundle. 

Among several interesting facts regarding holonomic submanifolds, we
recall that they are a natural play-ground for the Ribaucour transformation 
\cite{dt}. As an application of the so called vectorial Ribaucour 
transformation as given in \cite{dft2}, one can locally parametrically 
generate any proper holonomic submanifold in terms of a set of smooth 
functions whose Hessians are all diagonal with respect to the coordinate 
vector fields of a given orthogonal system of coordinates; see \cite{dft1} 
for details.

Since holonomic submanifolds have flat normal bundle, it is a standard
fact (see \cite{re}) that at each point $x\in M^n$ there exist a set 
of unique pairwise distinct normal vectors  
$\eta_i\in N_fM(x),\ 1\leq i\leq s=s(x)$,
and an associated orthogonal splitting of the tangent space as 
$$
T_xM=E_1(x)\oplus\cdots\oplus E_{s}(x)
$$
where
$$
E_i(x)=\{X\in T_xM:\alpha(X,Y)=\<X,Y\>\eta_i
\;\;\mbox{for all}\;\; Y\in T_xM\}.
$$
Hence, the second fundamental form of $f$ acquires the form
\be\label{sff}
\alpha(X,Y)(x)=\sum_{i=1}^s\< X^i,Y^i\> \eta_i
\ee
where $X\mapsto X^i$ is the orthogonal projection from $T_xM$ onto $E_i(x)$.

A submanifold $f\colon M^n\to\Q_c^N$ with flat normal bundle is said to be 
\emph{proper} if 
$s(x)=k$ is constant on $M^n$. If this is the case, then the maps 
$x\in M^n\mapsto\eta_i(x)$,  $1\leq i\leq k$, are smooth vector fields, 
called the \emph{principal normal vector fields} of $f$, and the 
distributions $x\in M^n\mapsto E_i(x)$, $1\leq i\leq k$, are also smooth. 

There are several conditions that imply that a submanifold of a space 
form has to be locally holonomic.  By locally we mean along connected
components of an open dense subset of the manifold. For instance, this 
is the case of any isometric immersion $f\colon M_c^n\to\Q_c^N$ with flat 
normal bundle of a manifold with the same constant sectional curvature 
as the ambient space form provided that index of relative nullity vanishes 
at any point; see \cite{dt} for a more general result.  Recall that the 
\emph{index of relative nullity} $\nu(x)$ of $f\colon M^n\to\Q_c^N$  
at $x\in M^n$ is the dimension of the \emph{relative nullity subspace} 
$\Delta(x)\subset T_xM$ given by
$$
\Delta(x)=\{X\in T_xM: \alpha(X,Y)=0\;\;\mbox{for all}\;\;Y\in T_xM\}.
$$

Isometric immersions $f\colon M_c^n\to\Q_{\tilde c}^{n+p}$
with sectional curvatures $c<\tilde{c}$ and in the least possible codimension 
$p=n-1$ have flat normal bundle and thus are always locally holonomic. 
This was already known to Cartan \cite{ca} who made an exhaustive study 
of the subject and  determined the degree of generality of such  submanifolds.
Moreover, being holonomic is also necessarily the case for $c>\tilde{c}$ 
but now under the  extra condition that the submanifold is
free of weak-umbilic points; see \cite{mo}.
\vspace{1ex}

In this note, we show that results discussed above still hold for 
isometric immersions of the larger class of Einstein manifolds. In fact, 
this turns out to be the case even in the presence of a constant positive  
index of relative nullity, thus in the case of submanifolds of manifolds 
with the same constant sectional curvature  the restriction mentioned above 
can be removed. 

\begin{theorem}\po\label{main}
Any isometric immersion $f\colon M^n\to\Q_c^N$ with flat normal bundle 
and proper of an Einstein manifold is locally holonomic. 
\end{theorem}

\proof It holds that $M^n$ is Einstein with 
$\Ric_M=\lambda\id$ if and only if the vector fields 
$$
\hat\eta_i=\eta_i-\frac{n}{2}H,\;\; 1\leq i\leq k,
$$ 
satisfy 
\be\label{hat}
\|\hat\eta_i\|^2= \frac{n^2}{4}\|H\|^2+c(n-1)-\lambda
\ee 
where $H$ denotes the mean curvature vector field of $f$. To see this, 
we first observe that for any isometric immersion $f\colon M^n\to\Q_c^N$ 
a straightforward computation of the  Ricci tensor using the Gauss 
equation yields 
$$
\Ric_M(X,Y)=c(n-1)\<X,Y\>+n\<\alpha(X,Y),H\>
-\sum_{j=1}^n\<\alpha(X,X_j),\alpha(Y,X_j)\>
$$
where $X_1,\ldots,X_n$ is an orthonormal tangent base.  
It follows easily using \eqref{sff} that
\be\label{r}
\Ric_M(X,Y)=c(n-1)\<X,Y\>+n\<\alpha(X,Y),H\>
-\sum_{i=1}^s\< X^i,Y^i\>\|\eta_i\|^2
\ee
for all $X,Y\in TM$. From \eqref{sff} and \eqref{r} we have that 
$\Ric_M=\lambda\id$ is equivalent to
$$
c(n-1)-\lambda=\|\eta_i\|^2-n\<H,\eta_i\>,\;\;1\leq i\leq k,
$$
and this in turn is equivalent to \eqref{hat}.  

We claim that at any point the vectors $\eta_i-\eta_j$ and 
$\eta_i-\eta_\ell$ are linearly independent 
if $i\neq j\neq\ell\neq i$. Assume to the contrary that 
$$
\eta_i-\eta_j=\mu(\eta_i-\eta_\ell)
$$
for some $\mu\neq 0$.  Then
$$
(1-\mu)\hat\eta_i=\hat\eta_j-\mu\hat\eta_\ell
$$
yields
$$
\|\hat\eta_i\|^2-2\mu\|\hat\eta_i\|^2+\mu^2\|\hat\eta_i\|^2=
\|\hat\eta_j\|^2-2\mu\<\hat\eta_j,\hat\eta_\ell\>+\mu^2\|\hat\eta_\ell \|^2.
$$
We have from \eqref{hat} that the $\hat\eta_j$'s are of equal length. Thus
$$
\<\hat\eta_j,\hat\eta_\ell\>=\|\hat\eta_j\|\|\hat\eta_\ell\|
$$
and hence $\hat\eta_j=\hat\eta_\ell$. This is a contradiction that proves 
the claim.

In order to conclude holonomicity it is a standard fact that it suffices 
to show is that the distributions $E_j^\perp$ are integrable for 
$1\leq j\leq k$. 
The Codazzi equation is easily seen to yield
\be\label{c1}
\<\nabla_XY,Z\>(\eta_i-\eta_j)=\<X,Y\>\nabla_Z^\perp\eta_i
\ee
and
\be\label{c2}
\<\nabla_XV,Z\>(\eta_j-\eta_\ell)=\<\nabla_VX,Z\>(\eta_j-\eta_i)
\ee
for all $X,Y\in E_i,Z\in E_j$ and $V\in E_\ell$ where
$1\leq i\neq j\neq\ell\neq i\leq k$.

It follows from \eqref{c1} that the $E_i$'s are integrable. 
Thus, it is sufficient to argue for the case $k\geq 3$.
In fact, it suffices to show that if $X\in E_i$ and $Y\in E_j$ then 
$[X,Y]\in E_\ell^\perp$ if $i\neq j\neq\ell\neq i$. 
We have from \eqref{c2} that
$$
\<\nabla_XY,Z\>(\eta_\ell-\eta_j)=\<\nabla_YX,Z\>(\eta_\ell-\eta_i)
$$
for any $Z\in E_\ell$. We obtain using the result of the claim that
$$
\<\nabla_X Y,Z\>=\<\nabla_Y X,Z\>=0,
$$ 
and this completes the proof.\vspace{2ex}\qed

Let $g\colon L^{n-s}\to\Q_c^N$, $1\leq s\leq n-1$, be an isometric 
immersion carrying a parallel flat normal subbundle 
$\mathcal{L}\subset N_gL$ of rank $s$. The \emph{generalized cylinder} 
determined by the subbundle $\pi:\mathcal{L}\to L^{n-s}$
is the $n$-dimensional submanifold $f\colon M^n\to\Q_c^N$ parametrized 
(at regular points) by means of the exponential map of $\Q_c^N$ as
$$
\gamma\in\mathcal{L}\mapsto\exp_{g(\pi(\gamma))}\gamma.
$$

We have that $\gamma\in\mathcal{L}$ is a regular point if and only if
$P=I-A^g_\gamma$ is nonsingular where $A^g_\gamma$ stands for the shape 
operator of $g$ corresponding to $\gamma$.  Also $N_fM=\mathcal{L}^\perp$, 
up to parallel identification along the fibers
of $\mathcal{L}$ that are contained in the relative nullity subspaces of $f$.  
Moreover, the relation between the second fundamental
forms of $f$ and $g$ is given by
$$
\alpha_f(X,Y)=\left(\alpha_g(X,PY)\right)_{\mathcal{L}^\perp}
$$
for all $X,Y\in TL$. It follows that $f$ has flat normal bundle if 
and only if $g$ has flat normal bundle.
\vspace{1ex}

The following result obtained in \cite{dft1} asserts that any submanifold 
with a relative nullity distribution $x\in M^n\mapsto\Delta(x)$  of constant 
dimension whose \emph{conullity} distribution $x\in M^n\mapsto\Delta^\perp(x)$
is integrable has to be a generalized cylinder.

\begin{proposition}\po\label{nul}  
Let $g\colon L^{n-s}\to\Q_c^N$ be an isometric immersion carrying a 
parallel flat normal subbundle $\mathcal{L}\subset N_gL$ of rank $s$  
such that 
$$
\{Y\in T_xL:(\alpha_g(Y,X))_{\mathcal{L}^\perp}=0\;\;\mbox{for all}\;\;X\in T_xL\}=0
$$
for any point $x\in M^n$.
Then the generalized cylinder over $g$ determined by  $\mathcal{L}$ has 
relative nullity of constant dimension $s$ and integrable conullity.

Conversely, any submanifold $f\colon M^n\to \Q_c^N$ with relative nullity 
distribution $\Delta$ of constant dimension $s$ and integrable conullity
arises this way locally. This means that $\mathcal{L}=\Delta|_L$ is a 
parallel flat normal subbundle of $g=f|_L$ for any integral leaf $L^{n-s}$ of 
the conullity and $f$ is locally an open neighborhood of $g(L)$ in the generalized 
cylinder over $g$ determined by $\mathcal{L}$.
\end{proposition}

We have the following consequence of Theorem \ref{main}.

\begin{corollary}\po\label{cor}
Let $f\colon M^n\to\Q_c^N$ be an isometric immersion with flat normal bundle 
of an Einstein manifold that is proper and has constant index of relative 
nullity $\nu=s\geq 1$.  Then $f$ is locally a generalized cylinder 
over a submanifold $g\colon L^{n-s}\to\Q_c^N$ with flat normal bundle.
\end{corollary}

\proof Notice that if the index of relative nullity is $\nu\geq 1$ at any 
point and $\nu>1$ at some point then it has to be constant since $f$ is proper.
The proof follows easily from Theorem~\ref{main} and Proposition~\ref{nul}.
\qed

\noindent Marcos Dajczer\\
IMPA -- Estrada Dona Castorina, 110\\
22460--320, Rio de Janeiro -- Brazil\\
e-mail: marcos@impa.br
\bigskip

\noindent Christos-Raent Onti\\
University of Ioannina \\
Department of Mathematics\\
Ioannina--Greece\\
e-mail: chonti@cc.uoi.gr
\bigskip

\noindent Theodoros Vlachos\\
University of Ioannina \\
Department of Mathematics\\
Ioannina--Greece\\
e-mail: tvlachos@uoi.gr

\end{document}